\documentclass[graybox]{svmult}

\usepackage{mathptmx}       
\usepackage{helvet}         
\usepackage{courier}        
\usepackage{type1cm}        
\usepackage{makeidx}         
\usepackage{graphicx}        
\usepackage{multicol}        
\usepackage[bottom]{footmisc}
\usepackage{amsmath, amssymb, multirow, booktabs, siunitx}

\makeindex             

\begin{document}
	
	\title*{Convergence Rates of a Fully Discrete Galerkin Scheme for the Benjamin--Ono Equation}
	\titlerunning{Convergence Rates of a Fully Discrete Galerkin Scheme for the BO Equation}
	\author{Sondre Tesdal Galtung}
	\institute{Sondre Tesdal Galtung \at Department of Mathematical Sciences, NTNU Norwegian University of Science and Technology \\ NO-7491 Trondheim, Norway \\ \email{sondre.galtung@ntnu.no}}
	%
	%
	\maketitle
	
	\abstract*{}
	
	\abstract{We consider a recently proposed fully discrete Galerkin scheme for the Benjamin--Ono equation which has been found to be locally convergent in finite time for initial data in $L^2(\mathbb{R})$. By assuming that the initial data is sufficiently regular we obtain theoretical convergence rates for the scheme both in the full line and periodic versions of the associated initial value problem. These rates are illustrated with some numerical examples.
	}
	
	\section{Background}
	We will in the following consider the Benjamin--Ono (BO) equation \cite{benjamin, ono} which serves as a generic model for weakly nonlinear long waves with nonlocal dispersion.
	Its initial value problem reads
	\begin{equation}
		\begin{cases}
		u_{t} + uu_{x} - \mathrm{H}u_{xx} = 0, & (t,x) \in (0,T] \times \mathbb{R}, \\
		u(0,x) = u_{0}(x), & x \in \mathbb{R},
		\end{cases}
		\label{BO}
	\end{equation}
	where $\mathrm{H}$ denotes the Hilbert transform defined by
	\begin{equation*}
	\mathrm{H}u(\cdot,x) := \textrm{p.v.} \frac{1}{\pi}\int_{\mathbb{R}} \frac{u(\cdot,x-y)}{y}\D y,
	\label{Hilbert}
	\end{equation*}
	for which p.v. denotes the Cauchy principal value.
	We may also consider the $2L$-periodic IVP for the BO equation
	\begin{equation}
	\begin{cases}
	u_{t} + uu_{x} - \mathrm{H}_{\text{per}}u_{xx} = 0, & (t,x) \in (0,T] \times \mathbb{T}, \\
	u(0,x) = u_{0}(x), & x \in \mathbb{T},
	\end{cases}
	\label{BO_per}
	\end{equation}
	where $\mathbb{T} := \mathbb{R}/2L\mathbb{Z}$, and $\mathrm{H}_{\text{per}}$ denotes the $2L$-periodic Hilbert transform defined by
	\begin{equation*}
	\mathrm{H}_{\text{per}}u(\cdot,x) := \textrm{p.v.} \frac{1}{2L} \int_{-L}^{L} u(\cdot,x-y) \cot\left(\frac{\pi}{2L}y\right) \D y.
	\label{Hilbert_per}
	\end{equation*}
	
	Based on a method for the Korteweg--de Vries equation due to Dutta and Risebro \cite{risebro2015}, Galtung \cite{galtung} proposed a fully discrete Crank--Nicolson Galerkin scheme for (\ref{BO}) where an inherent smoothing effect is used to prove convergence locally for initial data $u_{0}$ in $L^{2}(\mathbb{R})$ and a finite time $T$ which depends on $\|u_{0}\|_{L^2}$.
	
	The scheme for (\ref{BO}) is defined in the following way.
	First one discretizes a subset of the real line by dividing it in intervals of equal length $\Delta x$, $I_{j} = [x_{j-1}, x_{j}]$, where $x_{j} := j \Delta x,\, j \in \mathbb{Z}$.
	For the temporal discretization one analogously has $t_{n} = n \Delta t,\, n \in \{0,1,\dots,N\}$, for a discretization parameter $\Delta t$ such that $T = (N+1/2)\Delta t$.
	Let us also for convenience define $t_{n+1/2} := (t_{n} + t_{n+1})/2$.
	Consider now the following finite-dimensional subspace of the Sobolev space $H^{2}(\mathbb{R})$,
	\begin{equation}
	S_{\Delta x} = \{v \in H^{2}(\mathbb{R}) \:|\: v \in \mathbb{P}_{r}(I_{j}), j \in \mathbb{Z} \},
	\label{subspace}
	\end{equation}
	where $r \ge 2$ is a fixed integer and $\mathbb{P}_{r}(I)$ denotes the space of polynomials on the interval $I$ of degree less than or equal to $r$.
	Given $R > 0$ we define $\varphi \in C^{\infty}(\mathbb{R})$, for which the derivative is a cut-off function, satisfying the following conditions,
	\begin{enumerate}
		\item $1 \le \varphi(x) \le 2 + 2R$,
		\item $\varphi'(x) = 1$ for $|x| < R$,
		\item $\varphi'(x) = 0$ for $|x| \ge R+1$, and
		\item $0 \le \varphi'(x) \le 1$ for all $x$.
	\end{enumerate}
	This function plays a key role in establishing the previously mentioned smoothing effect for the scheme, and it may be chosen to be point-symmetric in $(0,\varphi(0))$.
	
	We need a reasonable approximation of $u_{0}$ in (\ref{BO}) as initial data $u^{0}$ for our scheme, and so we set $u^{0} = \mathrm{P}u_{0}$, where $\mathrm{P}$ is the $L^{2}$-projection on $S_{\Delta x}$.
	Now we define a sequence of approximations $\{u^{n}\}_{n=0}^{N}$ of the exact solution at each $t_{n}$ by the following procedure:
	find $u^{n+1} \in S_{\Delta x}$ such that
	\begin{equation}
	\left< u^{n+1}, \varphi v \right> - \frac{\Delta t}{2} \left< \left(u^{n+1/2}\right)^{2}, \left(\varphi v\right)_{x} \right> + \Delta t \left< \mathrm{H}\left(u^{n+1/2}\right)_{x}, \left(\varphi v\right)_{x} \right> = \left< u^{n}, \varphi v \right>,
	\label{scheme}
	\end{equation}
	for all $v \in S_{\Delta x}$, where $u^{0}$ is defined as before and $u^{n+1/2} := (u^{n} + u^{n+1})/2$.
	Here $\left<\cdot,\cdot\right>$ is the standard $L^{2}$-inner product.
	Note that the inner product $\left<\cdot,\cdot\varphi\right> =:\left<\cdot,\cdot\right>_{\varphi}$ defines a norm which we denote $\|\cdot\|_{2,\varphi}$.
	The nonlinearity appearing in the above implicit scheme calls for some form of iterative method to solve (\ref{scheme}) for each time step, and in \cite{galtung} the following linearized scheme is used,
	\begin{equation}
		\begin{cases}
		\left< w^{\ell+1}, \varphi v \right> - \frac{\Delta t}{2} \left< \left(\frac{w^{\ell}+u^{n}}{2}\right)^{2}, (\varphi v)_{x} \right> + \Delta t \left< \left(\mathrm{H} \frac{w^{\ell+1}+u^{n}}{2}\right)_{x}, (\varphi v)_{x} \right> = \left< w^{\ell}, \varphi v \right>, \\
		w^{0} = u^{n},
		\end{cases}
		\label{iteration}
	\end{equation}
	which is to hold for all $v \in S_{\Delta x}$.
	By assuming a CFL condition of the type $\Delta t = \mathrm{O}(\Delta x^{2})$ the above iteration is shown to converge to the solution $u^{n+1}$ of (\ref{scheme}).
	From this one can show that there exists $T > 0$ such that $u^{\Delta x}$, which is a piecewise linear interpolation of each $u^{n}$, belongs to the space $L^{2}(0,T;H_{\mathrm{loc}}^{1/2}(\mathbb{R}))$.
	Then compactness arguments yield the convergence result.
	
	Because a monotone increasing cut-off function is incompatible with the periodicity of (\ref{BO_per}) one cannot use the same arguments to prove convergence for $L^{2}$-initial data in this case, and so other tools are called for when considering low regularity initial data for the periodic BO equation.
	However, in this study  we will assume the initial data to be as regular as needed, and so we will consider the convergence rate of the method in best-case scenarios.
	The established well-posedness of the BO equation for these more regular spaces then guarantees that the exact solution at all times is at least as regular as the initial data.
	This will even make us able to consider the periodic IVP (\ref{BO_per}) using a slightly adapted scheme where we have simply replaced the cut-off function $\varphi$ with $1$ wherever it appears.
	
	In the upcoming analysis we need some preliminary estimates for polynomial approximations in finite element spaces.
	For a function $v \in S_{\Delta x}$ we have the following inverse inequalities
	\begin{align}
	|v|_{W^{k,\infty}(\mathbb{R})} &\le \frac{C}{(\Delta x)^{1/2}} |v|_{H^{k}(\mathbb{R})}, \qquad k = 0, 1, \label{inverse_1}\\
	|v|_{H^{k+1}(\mathbb{R})} &\le \frac{C}{\Delta x} |v|_{H^{k}(\mathbb{R})}, \qquad k = 0, 1, \label{inverse_2}
	\end{align}
	where the constant $C$ is independent of $v$ and $\Delta x$.
	Both here and in the following, $|\cdot|_{W^{k,p}(\mathbb{R})}$ denotes the seminorm of the Sobolev space $W^{k,p}(\mathbb{R})$ for which $H^{k}(\mathbb{R}) := W^{k,2}(\mathbb{R})$.
	The reader is referred to \cite[p. 142]{ciarlet} for a proof of the above inequalities.
	
	Let us now consider two projections $\mathrm{P} : L^{2}(\mathbb{R}) \to S_{\Delta x}$ and $\mathrm{P}_{\varphi} : L^{2}(\mathbb{R}) \to S_{\Delta x}$ defined respectively by
	\begin{equation}
		\int_{\mathbb{R}} \left(\mathrm{P}u - u\right) v \,\D x = 0,\quad v \in S_{\Delta x}, \label{proj}
	\end{equation}
	and
	\begin{equation}
		\int_{\mathbb{R}} \left(\mathrm{P}_{\varphi}u - u\right) \varphi v \,\D x = 0,\quad v \in S_{\Delta x}. \label{proj_phi}
	\end{equation}
	For these projections applied to a function $u \in H^{2}(\mathbb{R})$, we have the bounds
	\begin{align}
		\begin{split}
		\|\mathrm{P}_{0}u\|_{L^2(\mathbb{R})} &\le C \|u\|_{L^{2}(\mathbb{R})}, \\
		\|\mathrm{P}_{0}u\|_{H^1(\mathbb{R})} &\le C \|u\|_{H^{1}(\mathbb{R})}, \\
		\|\mathrm{P}_{0}u\|_{H^2(\mathbb{R})} &\le C \|u\|_{H^{2}(\mathbb{R})},		
		\end{split}
		\label{proj_bounds}
	\end{align}
	where $\mathrm{P}_{0}$ denotes either of the two projections and $C$ is a constant which is independent of $\Delta x$.
	These bounds can be derived from the norm equivalence in a finite-dimensional space and the definitions of these projections.
	
	We also have the following polynomial approximation error estimate on the discretized domain $\Omega$. Given $u \in H^{l+1}(\Omega)$, $0 \le m \le l$ and $s := \min\{l,r\}$, then
	\begin{align}
		|\mathrm{P}_{0}u - u|_{H^{m}(\Omega)} \le C \Delta x^{s+1-m} |u|_{ H^{s+1}(\Omega) }, \quad m = 0,1,2, \label{proj_err}
	\end{align}
	where again $\mathrm{P}_{0}$ denotes either of the two projections and $C$ is a constant not depending on $\Delta x$.
	For a proof of \eqref{proj_err} for $\mathrm{P}$ we refer to \cite[p. 98]{quarteroni} and the result for $\mathrm{P}_{\varphi}$ follows from an adaption of the same proof.
	
	The following properties of the Hilbert transform, which can be found in \cite[p. 317]{grafakos}, are also useful:
	\begin{align*}
	\left< \mathrm{H}u, v\right> &= -\left< u, \mathrm{H}v\right> \text{ for } u, v \in L^2(\mathbb{R}), \\
	(\mathrm{H}u)_{x} &= \mathrm{H}u_{x}, \\
	\|\mathrm{H}u\|_{L^2(\mathbb{R})} &= \|u\|_{L^2(\mathbb{R})}.
	\end{align*}
	Note that these properties hold analogously for the $2L$-periodic Hilbert transform $\mathrm{H}_{\text{per}}$ on $\mathbb{T}$ with $L^2(\mathbb{T}) = L^2([-L,L])$, except that $\|\mathrm{H}_{\text{per}}u\|_{L^2(\mathbb{T})} \le \|u\|_{L^2(\mathbb{T})}$.
	
	\section{Analysis of convergence rates}
	
	In the following we want to consider the $L^2$-norm of the difference $u^{n}-u(t_{n})$, and we will do so by decomposing the error as
	\begin{equation*}
	u^{n}-u(t_{n}) = (u^{n} - \mathrm{P}_{0}u(t_{n})) + (\mathrm{P}_{0}u(t_{n}) - u(t_{n})) =: \tau^{n} + \rho^{n},
	\end{equation*}
	and we will use the notation $w^{n} := \mathrm{P}_{0}u(t_{n})$ for the sake of brevity.
	Here $\mathrm{P}_{0} = \mathrm{P}_{\varphi}$ in the full line case, and $\mathrm{P}_{0} = \mathrm{P}$ for the periodic case.
	For $\rho^{n}$ we already have estimates for the $L^2$-norm by virtue of \eqref{proj_err}, and so it remains to estimate the norm of $\tau^{n}$.
	As the analysis is similar for the full line and periodic problems, we will give detailed estimates for the former case and only indicate the main differences between the two for the latter case.
	Note that in the following, $C$ will denote a constant which exact value is of no importance.
	Similarly, $C(R)$ will denote such a constant which depends on $R$ and so on.
	When we write  e.g., $L^2$ it is understood from context if we are referring to $L^2(\mathbb{R})$ or $L^2(\mathbb{T}) = L^2([-L,L])$.
	For both the full line and periodic case we have the following result which is proved in the next subsections.
	\begin{theorem}
		Given sufficiently regular initial data $u_{0}$, say $u_{0} \in H^{\max\{r+1,6\}}$, for the IVP of the BO equation, we have the following convergence rate for the fully discrete Galerkin scheme described in the previous section,
		\begin{equation}
			\|u^{n} - u(t_{n})\|_{L^2} = O(\Delta x^{r-1} + \Delta t^2), \quad n = 0,\dots,N.
			\label{rates}
		\end{equation}
		\label{main thm}
	\end{theorem}
	
	\subsection{Full line problem}
	From multiplying \eqref{BO} by $\varphi v$, where $v \in H^{2}$, and integrating by parts we get
	\begin{equation}
		\left<u_{t}(t), \varphi v\right> -\frac{1}{2}\left<u(t)^2,(\varphi v)_{x} \right> + \left< \mathrm{H}u_{x}(t), (\varphi v)_{x} \right> = 0, \quad t \in (0,T].
		\label{weak_form}
	\end{equation}
	From \eqref{scheme}, \eqref{weak_form} and \eqref{proj_phi} we are able to write
	\begin{align*}
		\left<\frac{\tau^{n+1}-\tau^{n}}{\Delta t}, \varphi v\right> &= \left<\frac{u^{n+1}-u^{n}}{\Delta t}, \varphi v\right> - \left<\frac{w^{n+1}-w^{n}}{\Delta t}, \varphi v\right> \\
		&= \left<\frac{u^{n+1}-u^{n}}{\Delta t}, \varphi v\right> - \left<u_{t}(t_{n+1/2}), \varphi v\right> \\
		&\quad+ \left<\underbrace{ u_{t}(t_{n+1/2}) - \frac{u(t_{n+1})-u(t_{n})}{\Delta t}}_{\kappa^{n+1/2}}, \varphi v\right> \\
		&= -\frac{1}{2}\left<(u^{n+1/2})^2-u(t_{n+1/2})^2,(\varphi v)_{x} \right> \\
		&\quad+ \left< \mathrm{H}(u^{n+1/2}_{x}-u_{x}(t_{n+1/2})), (\varphi v)_{x} \right> + \left< \kappa^{n+1/2},\varphi v\right>,
	\end{align*}
	for $v \in S_{\Delta x}$.
	As we are now considering $u$ evaluated at $t_{n+1/2}$ we cannot use the previous decomposition of the error directly, but we instead write
	\begin{equation*}
		u^{n+1/2} - u(t_{n+1/2}) = \tau^{n+1/2} + \rho^{n+1/2} + \underbrace{\frac{u(t_{n+1})+ u(t_{n})}{2} - u(t_{n+1/2})}_{\sigma^{n+1/2}}.
	\end{equation*}
	Then we may rewrite part of the nonlinear term as
	\begin{align*}
		(u^{n+1/2})^2 - u(t_{n+1/2})^2 &= (\tau^{n+1/2})^2 + 2\tau^{n+1/2}w^{n+1/2} + (w^{n+1/2})^2 - u(t_{n+1/2})^2 \\
		&= (\tau^{n+1/2})^2 + 2\tau^{n+1/2}w^{n+1/2} \\
		&\quad+ (w^{n+1/2} + u(t_{n+1/2}))(\rho^{n+1/2}+\sigma^{n+1/2}).
	\end{align*}
	In the following we want to use $\tau^{n+1/2} \in S_{\Delta x}$ as test function,
	and from integrating by parts we get the following relevant identities,
	\begin{align*}
		\left<(\tau^{n+1/2})^2,(\varphi \tau^{n+1/2})_{x}\right> &= -\frac{1}{3}\left<(\tau^{n+1/2})^3,\varphi_{x}\right>, \\
		2\left<\tau^{n+1/2}w^{n+1/2},(\varphi \tau^{n+1/2})_{x} \right> &= -\left<(\tau^{n+1/2})^2,\varphi w^{n+1/2}_{x}\right> + \left<(\tau^{n+1/2})^2,\varphi_{x} w^{n+1/2}\right>.
	\end{align*}
	
	Inserting this in the previous equations we get
	\begin{align*}
		\frac{1}{2}\|\tau^{n+1}\|^{2}_{2,\varphi} &= \frac{1}{2}\|\tau^{n}\|^{2}_{2,\varphi} + \Delta t \left[ -\frac{1}{6} \left< (\tau^{n+1/2})^3,\varphi_{x} \right> -\frac{1}{2} \left< (\tau^{n+1/2})^2,\varphi w^{n+1/2}_{x} \right> \right. \\
		&\quad+ \left. \frac{1}{2} \left< (\tau^{n+1/2})^2,\varphi_{x} w^{n+1/2} \right> \right. \\
		&\quad+ \left. \frac{1}{2} \left<(w^{n+1/2}+u(t_{n+1/2}))(\rho^{n+1/2}+\sigma^{n+1/2}),(\varphi \tau^{n+1/2})_{x} \right> \right. \\
		&\quad- \left. \left<\mathrm{H}\tau^{n+1/2}_{x}, (\varphi \tau^{n+1/2})_{x}\right>  -\left<\mathrm{H}\rho^{n+1/2}_{x}, (\varphi \tau^{n+1/2})_{x}\right> \right. \\
		&\quad- \left. \left<\mathrm{H}\sigma^{n+1/2}_{x}, (\varphi \tau^{n+1/2})_{x}\right> + \left< \kappa^{n+1/2}, \varphi \tau^{n+1/2} \right> \vphantom{\frac{1}{6}} \right].
	\end{align*}
	From the commutator estimates presented in \cite{galtung} we have the inequalities
	\begin{equation*}
		\left< \mathrm{H}w_{x}, (\varphi w)_{x} \right> \ge \left\|\sqrt{\varphi_{x}} D^{1/2} w \right\|^{2}_{L^2} - \widetilde{C} \|w\|_{L^2}^2,
	\end{equation*}
	and
	\begin{equation*}
		\left<w^3,\varphi_{x}\right> \le \left\|\sqrt{\varphi_{x}} D^{1/2} w \right\|^{2}_{L^2} + C(1+\|w\|^{2}_{L^2})\|w\|^{2}_{L^2}
	\end{equation*}
	for $w \in H^{2}$.
	By inserting these in the preceding identity and using the $L^2$-isometry of the Hilbert transform we obtain
	\begin{align*}
		\frac{1}{2}\|\tau^{n+1}\|^{2}_{2,\varphi} \hspace{-1em}&\hspace{1em}+ \Delta t \left\|\sqrt{\varphi_{x}} D^{1/2} \tau^{n+1/2} \right\|^{2}_{L^2} - \Delta t \widetilde{C} \|\tau^{n+1/2}\|_{2,\varphi}^2 \\
		&\le \frac{1}{2}\|\tau^{n}\|^{2}_{2,\varphi} + \frac{\Delta t}{3} \left\|\sqrt{\varphi_{x}} D^{1/2} \tau^{n+1/2} \right\|^{2}_{L^2} \\
		&\quad+ \Delta t \left[C(1+\|u^{n+1/2}-w^{n+1/2}\|^{2}_{L^2})\|\tau^{n+1/2}\|_{2,\varphi}^2 \right. \\
		&\quad+ \left. \frac{1}{2}\|w^{n+1/2}_{x}\|_{L^{\infty}} \|\tau^{n+1/2}\|_{2,\varphi}^{2} + \frac{1}{2}\|w^{n+1/2}\|_{L^{\infty}} \|\tau^{n+1/2}\|_{2,\varphi}^{2} \right. \\
		&\quad+ \left. \frac{1}{2}\|w^{n+1/2}_{x} + u_{x}(t_{n+1/2})\|_{L^{\infty}}(\|\rho^{n+1/2}\|_{2,\varphi} + \|\sigma^{n+1/2}\|_{2,\varphi})\|\tau^{n+1/2}\|_{2,\varphi} \right. \\
		&\quad+ \left. \frac{1}{2}\|w^{n+1/2} + u(t_{n+1/2})\|_{L^{\infty}}(\|\rho_{x}^{n+1/2}\|_{2,\varphi} + \|\sigma_{x}^{n+1/2}\|_{2,\varphi})\|\tau^{n+1/2}\|_{2,\varphi} \right. \\
		&\quad+ \left. C_{R} \|\rho^{n+1/2}_{x}\|_{L^{2}} (\|\tau^{n+1/2}\|_{2,\varphi} + C_{R} \|\tau^{n+1/2}_{x}\|_{L^{2}}) \right. \\
		&\quad+ \left. C_{R} \|\sigma^{n+1/2}_{xx}\|_{L^2}\|\tau^{n+1/2}\|_{2,\varphi} + C_{R}\|\kappa^{n+1/2}\|_{L^2}\|\tau^{n+1/2}\|_{2,\varphi} \right].
	\end{align*}
	From the Sobolev inequality $\|w \|_{L^{\infty}(\mathbb{R})} \le \|w\|_{H^1(\mathbb{R})}$, the Cauchy--Schwarz inequality, \eqref{inverse_1} and reordering we then obtain
	\begin{align*}
	\frac{1}{2}\|\tau^{n+1}\|^{2}_{2,\varphi} \hspace{-1em}&\hspace{1em}+ \frac{2 \Delta t}{3} \left\|\sqrt{\varphi_{x}} D^{1/2} \tau^{n+1/2} \right\|^{2}_{L^2} \\
	&\le \frac{1}{2}\|\tau^{n}\|^{2}_{2,\varphi} + \Delta t C_{R} \left[(1+\|u^{n+1/2}\|^{2}_{L^2}+\|w^{n+1/2}\|^{2}_{L^2})\|\tau^{n+1/2}\|_{2,\varphi}^2 \right. \\
	&\quad+ \left. \|w^{n+1/2}\|_{H^{2}} \|\tau^{n+1/2}\|_{2,\varphi}^{2} + \frac{1}{2}\|w^{n+1/2}\|_{H^{1}} \|\tau^{n+1/2}\|_{2,\varphi}^{2} \right. \\
	&\quad+ \left. (\|w^{n+1/2}\|_{H^2} + \|u(t_{n+1/2})\|_{H^2})(\|\rho^{n+1/2}\|_{L^2} + \|\sigma^{n+1/2}\|_{L^2})\|\tau^{n+1/2}\|_{2,\varphi} \right. \\
	&\quad+ \left. (\|w^{n+1/2}\|_{H^1} + u(t_{n+1/2})\|_{H^1})(\|\rho_{x}^{n+1/2}\|_{L^2} + \|\sigma_{x}^{n+1/2}\|_{L^2})\|\tau^{n+1/2}\|_{2,\varphi} \right. \\
	&\quad+ \left. \|\rho^{n+1/2}_{x}\|_{L^{2}} (\|\tau^{n+1/2}\|_{2,\varphi} + \frac{1}{\Delta x} \|\tau^{n+1/2}\|_{2,\varphi}) \right. \\
	& \quad+ \left. \|\sigma^{n+1/2}_{xx}\|_{L^2}\|\tau^{n+1/2}\|_{2,\varphi} + \|\kappa^{n+1/2}\|_{L^2}\|\tau^{n+1/2}\|_{2,\varphi} \right].
	\end{align*}
	The following result is a part of Lemma 4.1 in \cite{galtung} and will be of use.
	\begin{lemma}
		Let $u^n$ be the solution of \eqref{scheme} and assume furthermore that the scheme fulfills a CFL condition of the form $\Delta t^2 / \Delta x^3 \le \tilde{C}$, where $\tilde{C}$ is a constant depending on $\|u_0\|_{L^2}$. 
		Then $\|u^{n}\|_{L^2} \le C(\|u_{0}\|_{L^{2}})$ for $n = 0,...,N$.
		\label{L2 bound}
	\end{lemma}
	Using Lemma \ref{L2 bound}, \eqref{proj_bounds}, Cauchy's inequality and dropping the second term on the left hand side we get
	\begin{align*}
	\|\tau^{n+1}\|^{2}_{2,\varphi} &\le \|\tau^{n}\|^{2}_{2,\varphi} + \Delta t C(u,R) \left[\|\tau^{n+1}\|_{2,\varphi}^2 + \|\tau^{n}\|_{2,\varphi}^2 + \|\rho^{n+1/2}\|_{L^2}^2 + |\rho^{n+1/2}|_{H^1}^2  \right. \\
	&\quad+ \left. \frac{1}{\Delta x^2}|\rho^{n+1/2}|_{H^1}^2 + \|\sigma^{n+1/2}\|_{H^2}^2 + \|\kappa^{n+1/2}\|_{L^2}^2 \right],
	\end{align*}
	which implies
	\begin{align*}
		(1-\Delta t C(u,R))\|\tau^{n+1}\|^{2}_{2,\varphi} \le (1+\Delta t C(u,R))\|\tau^{n}\|^{2}_{2,\varphi} + \Delta t C(u,R) S_{n},
	\end{align*}
	where we have the remainder term
	\begin{equation*}
		S_{n} = \|\rho^{n+1/2}\|_{L^2}^2 + |\rho^{n+1/2}|_{H^1}^2 + \frac{1}{\Delta x^2}|\rho^{n+1/2}|_{H^1}^2 + \|\sigma^{n+1/2}\|_{H^2}^2 + \|\kappa^{n+1/2}\|_{L^2}^2.
	\end{equation*}
	We will assume $\Delta t$ small enough that the left hand side of the previous inequality is strictly positive, say $1-\Delta t C(u,R)) \ge 1/2$.
	From Taylor's formula with integral remainder
	we can derive the following estimate for the seminorms of $\sigma^{n+1/2}$,
	\begin{equation}
	|\sigma^{n+1/2}|_{H^{k}}^{2} \le C \Delta t^3 \int_{t_{n}}^{t_{n+1}}|u_{tt}(s)|_{H^{k}}^2\,\D s,
	\label{sigma_est}
	\end{equation}
	and the $L^2$-norm of $\kappa^{n+1/2}$,
	\begin{equation}
	\|\kappa^{n+1/2}\|_{L^2}^2 \le C \Delta t^3 \int_{t_{n}}^{t_{n+1}}\|u_{ttt}(s)\|_{L^{2}}^2\,\D s.
	\label{kappa_est}
	\end{equation}
	Then we may estimate the remainder term using \eqref{proj_err}, \eqref{sigma_est} and \eqref{kappa_est},
	\begin{align*}
		S_{n} &\le C \Delta x^{2(r+1)}(|u(t_{n})|_{H^{r+1}} + |u(t_{n+1})|_{H^{r+1}}) + \frac{C\Delta x^{2r}}{\Delta x^2} (|u(t_{n})|_{H^{r+1}}^2 + |u(t_{n+1})|_{H^{r+1}}^2) \\
		&\quad+ C \Delta t^3 \int_{t_{n}}^{t_{n+1}} \|u_{tt}(s)\|_{H^2}^2\,\D s + C \Delta t^3 \int_{t_{n}}^{t_{n+1}} \|u_{ttt}(s)\|_{L^2}^2\,\D s \\
		&= C \Delta x^{2(r-1)} \sup_{0 \le t \le T} |u(t)|_{H^{r+1}}^2 + C \Delta t^3 \left( \int_{t_{n}}^{t_{n+1}} \|u_{tt}(s)\|_{H^2}^2\,\D s + \int_{t_{n}}^{t_{n+1}} \|u_{ttt}(s)\|_{L^2}^2\,\D s \right).
	\end{align*}
	This yields
	\begin{align*}
		\|\tau^{n}\|^{2}_{2,\varphi} &\le \left(\frac{1+C\Delta t}{1-C\Delta t}\right)^{n} \|\tau^{0}\|^{2}_{2,\varphi} + \Delta t C \sum_{j=0}^{n-1} \left(\frac{1+C\Delta t}{1-C\Delta t}\right)^{n-j} S_{j} \\
		&\le \E^{4C T} \|\tau^{0}\|^{2}_{2,\varphi} + \Delta t \E^{4C T} \sum_{j=0}^{n-1} S_{j} \\
		&\le T C(u,R,T) \Delta x^{2(r-1)} + C(T) \Delta t^4 \left( \int_{0}^{T} \|u_{tt}(s)\|_{H^2}^2\,\D s +  \int_{0}^{T} \|u_{ttt}(s)\|_{L^2}^2\,\D s \right) \\
		&= C(u,R,T)(\Delta x^{2(r-1)} + \Delta t^4).
	\end{align*}
	To ensure that the above norms are bounded we assume that $u_0 \in H^{s}(\mathbb{R})$, $s \ge \max\{r+1,6\}$, see Theorems 5.3.1 and 9.1 in \cite{abdelouhab}.
	Then we have
	\begin{equation*}
		\|\tau^{n}\|_{L^2} \le \|\tau^{n}\|_{2,\varphi} \le C(u,R,T) (\Delta x^{r-1} + \Delta t^2),
	\end{equation*}
	where we have employed \eqref{proj_err} to deduce
	\begin{equation*}
		\|\tau^0\|_{L^2} \le \|\mathrm{P}u_{0} - u_{0}\|_{L^2} + \|u_{0} - \mathrm{P}_{\varphi}u_{0}\|_{L^2} \le C \Delta x^{r+1},
	\end{equation*}
	and if one in the original scheme instead had set $u^0 := \mathrm{P}_{\varphi}u_{0}$, then one would have $\tau^0 = 0$ directly.
	From this and \eqref{proj_err} we get
	\begin{equation*}
		\|u^{n}-u(t_{n})\|_{L^2} \le \|\tau^{n}\|_{L^2} + \|\rho^{n}\|_{L^2}  \le C(u,R,T)(\Delta x^{r-1} + \Delta t^2), \quad n = 1,\dots,N,
	\end{equation*}
	which proves Theorem \ref{main thm} for the full line case.
	
	\subsection{Periodic problem}
	For the $2L$-periodic case we follow the steps made for the real line case, but without the cut-off function $\varphi$ involved in the scheme and all inner products now act on $[-L,L]$.
	In this case it is straightforward to check that the $L^2$-norm of the fully discrete solution $u^{n}$ is conserved, simply by choosing $v = u^{n+1/2}$ in the adapted version of \eqref{scheme}, integrating by parts and applying the skew-symmetry of the Hilbert transform and the periodicity of $u^n$.
	The existence and uniqueness for solutions the adapted version of the iterative scheme \eqref{iteration} can be done analogously to the original version.
	In this case we do not have the commutator estimates which were used to bound the terms $\left< \mathrm{H}\tau^{n}_{x}, (\varphi \tau^{n})_{x} \right>$ and $\left<(\tau^{n})^2, (\varphi \tau^{n})_{x}\right>$ by $\|\tau^{n}\|_{2,\varphi}^2$, but since these now appear as respectively $\left< \mathrm{H}_{\text{per}}\tau^{n}_{x}, \tau^{n}_{x} \right>$ and $\left<(\tau^{n})^2, \tau^{n}_{x}\right>$ we use the skew-symmetricity of $\mathrm{H}_{\text{per}}$ and the periodicity of $\tau^{n}$ to conclude that they both vanish. 
	Apart from this, one proceeds similarly to obtain the estimate \eqref{rates} for the periodic problem.
	Note that by obtaining this estimate we have proved the convergence of the scheme in the periodic case given sufficiently regular initial data using a stability and consistency argument.
	
	\section{Numerical experiments}
	In order to verify the convergence rates numerically we applied the fully discrete schemes to the problems \eqref{BO} and \eqref{BO_per}.
	Inspired by \cite{risebro2015} we define the subspace $S_{\Delta x}$ as follows.
	Let $f$ and $g$ be the functions
	\begin{align*}
	f(y) &= \begin{cases}
	1 + y^2 (2|y|-3), & |y| \le 1, \\
	0, & |y| > 1,
	\end{cases} \\
	g(y) &= \begin{cases}
	y(1-|y|)^2, & |y| \le 1, \\
	0, & |y| > 1.
	\end{cases}
	\end{align*}
	For $j \in \mathbb{Z}$ we define the basis functions
	\begin{align*}
	v_{2j}(x) = f\left(\frac{x-x_{j}}{\Delta x}\right), \quad
	v_{2j+1}(x) = g\left(\frac{x-x_{j}}{\Delta x}\right),
	\end{align*}
	where $x_{j} = j\Delta x$.
	Then $\{v_{j}\}_{-M}^{M}$ spans a $4M+2$ dimensional subspace of $H^{2}(\mathbb{R})$.
	In the following we define $N := 2M$, which is the number of elements used in the approximation.
	Note that for this choice we have $r=3$ in \eqref{subspace}, and so we expect convergence rates of order $O(\Delta x^2 + \Delta t^2)$.
	
	To approximate the full line for \eqref{BO} we have chosen to consider a finite interval with periodic boundary conditions, and we claim that this is a reasonable approximation as long as the approximate and exact solutions are close to zero at the endpoints, simulating the decay at infinity on the real line, which is the case for our examples.
	We have chosen to set $\Delta t = O(\Delta x)$, contrary to the assertion $\Delta t = O(\Delta x^{2})$ from the theory, as smaller time steps did not lead to significant improvement in the accuracy of  the approximations.
	In the iteration \eqref{iteration} to obtain $u^{n+1}$ we chose the stopping condition $\|w^{\ell+1}-w^{\ell}\|_{L^{2}} \le 0.002 \Delta x \|u^{n}\|_{L^{2}}$.
	The integrals involved in the Hilbert transforms were computed with seven and eight point Gauss--Legendre quadrature rules respectively for the inner Cauchy principal value integral and the outer integral appearing in the inner product.
	For $t = n\Delta t$, we set $u_{\Delta x}(x,t) = u^{n}(x,t) = \sum_{j=-M}^{M} u_{j}^{n} v_{j}(x)$.
	We have measured the relative error $E := \|u_{\Delta x} - u\|_{L^{2}}/\|u\|_{L^{2}}$ of the numerical approximation compared to the exact solution $u$, where the $L^{2}$-norms were computed with the trapezoidal rule in the grid points $x_{j}$ of the finest grid considered.
	
	\subsection{Full line problem}
	
\begin{figure}
		\centering
		\includegraphics[width=0.75\textwidth]{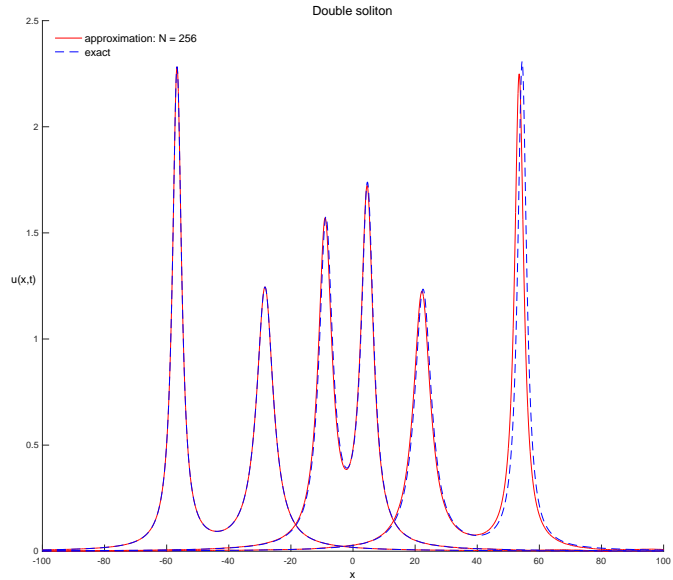}
		\caption{Numerical approximation for $N=256$ and exact solution for $t=0$, $90$ and $180$, respectively positioned from left to right in the plot, for full line problem with periodic boundary conditions. This figure is reproduced from \cite{galtung}}
		\label{twosolBO_per}
	\end{figure}

	A solution to this problem is the double soliton given by
	\begin{equation*}
	u_{s2}(x,t) = \frac{4c_{1}c_{2} \left( c_{1}\lambda_{1}^{2} + c_{2}\lambda_{2}^{2} + (c_{1}+c_{2})^{2}c_{1}^{-1}c_{2}^{-1}(c_{1}-c_{2})^{-2} \right)}{\left( c_{1}c_{2}\lambda_{1}\lambda_{2}-(c_{1}+c_{2})^{2}(c_{1}-c_{2})^{-2} \right)^{2} + \left( c_{1}\lambda_{1} + c_{2}\lambda_{2} \right)^{2}},
	\label{twoSol}
	\end{equation*}
	where $\lambda_{1} := x - c_{1}t - d_{1}$ and $\lambda_{2} := x - c_{2}t - d_{2}$. When $c_{2} > c_{1}$ and $d_{1} > d_{2}$, this equation represents a tall soliton overtaking a smaller one while moving to the right.
	We applied the fully discrete scheme with initial data $u_{0}(x) = u_{s2}(x,0)$ and parameters $c_{1} = 0.3$, $c_{2} = 0.6$, $d_{1} = -30$, and $d_{2} = -55$.
	The time step was set to $\Delta t = 0.5\Delta x/\|u_{0}\|_{L^{\infty}}$ and the numerical solutions were computed for $t = 90$ and $t = 180$, that is, during and after the taller soliton overtakes the smaller one.
	To approximate the full line problem we set the domain to $[-100,100]$ with the aforementioned periodic boundary condition, and based on this domain we chose the weight function $\varphi(x) = 120 + x$ for all experiments in this setting.
	The results are presented in Table \ref{twosol_norms_Hper} and a comparison between the approximation for $N = 256$ and the exact solution is shown in Fig. \ref{twosolBO_per}.
	These results for the full line problem are also presented as numerical examples for this scheme in \cite{galtung}.
	\begin{table}
		\caption{Relative $L^{2}$-error at $t=90$ and $t=180$ for full line problem with initial data $u_{s2}$ and periodic boundary conditions}
		\centering
		\begin{tabular}{ccccc}
			\toprule
			& \multicolumn{2}{c}{$t = 90$} & \multicolumn{2}{c}{$t = 180$} \\
			$N$ & $E$ & rate & $E$ & rate \\
			\midrule
			128 & 0.01844 & \multirow{2}{*}{-1.45} & 0.11959 & \multirow{2}{*}{-1.32} \\
			256 & 0.05021 & \multirow{2}{*}{1.58} & 0.29755 & \multirow{2}{*}{1.75} \\
			512 & 0.01678 & \multirow{2}{*}{0.68} & 0.08869 & \multirow{2}{*}{0.74} \\
			1024 & 0.01044 & \multirow{2}{*}{1.16} & 0.05295 & \multirow{2}{*}{2.35} \\
			2048 & 0.00467 & \multirow{2}{*}{0.08} & 0.01040 & \multirow{2}{*}{0.89} \\
			4096 & 0.00442 & & 0.00561 & \\
		\end{tabular}
		\label{twosol_norms_Hper}
	\end{table}
	The plot shows that the numerical approximation appears to be close to the exact solution and this is confirmed by the errors which are decreasing from $N = 256$ onwards, but not with a consistent rate.
	According to our analysis we should expect a convergence rate of 2, but at $t=180$ it varies from slightly below 1 to slightly above 2.
	As pointed out in \cite{galtung}, this is a complicated numerical example since one has to approximate the nonlinear interaction between two solitons.
	Moreover, approximating the full line by a periodic finite interval could also be contributing to the error and thus we are led to believe that the method applied to the periodic problem will yield results which are in better agreement with theory.
	
	\subsection{Periodic problem}
	
\begin{figure}
		\centering
		\includegraphics[width=0.75\textwidth]{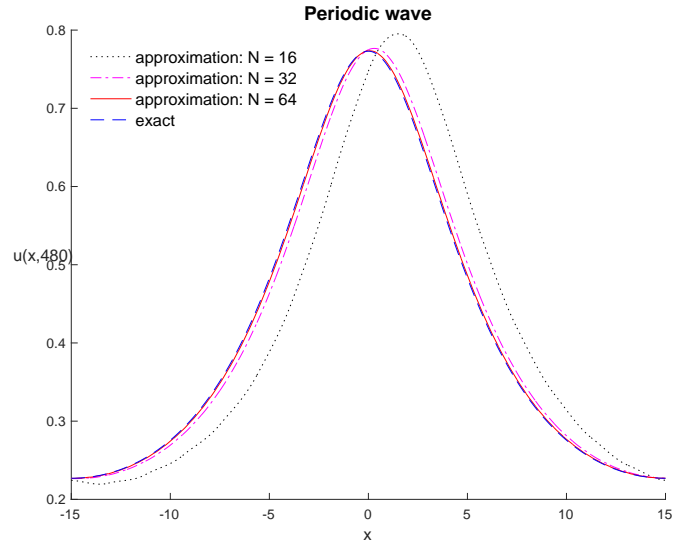}
		\caption{Exact and numerical solutions of the $2L$-periodic problem at $t=480$ for element numbers $N=16,32$ and $64$, with $L=15$ and initial data $u_{p1}$}
		\label{onewave_comp}
	\end{figure}

	In our second example we consider the Cauchy problem for the $2L$-periodic BO equation \eqref{BO_per}.
	In this case there exists a 2$L$-periodic single wave solution that tends to a single soliton as the period goes to infinity, given by
	\begin{equation*}
	u_{p1}(x,t) = \frac{2 c \delta}{1 - \sqrt{1 - \delta^{2}}\cos\left(c\delta(x-ct)\right)}, \quad \delta = \frac{\pi}{c L}.
	\label{oneWave}
	\end{equation*}
	We applied the scheme with initial data $u_{0}(x) = u_{p1}(x,0)$ with parameters $c = 0.25$ and $L = 15$.
	The time step was set to $\Delta t = 0.5 \Delta x$ and the approximate solution was computed for $t = 480$, which is four periods for the exact solution.
	As previously mentioned we do not have a weight function in this setting, which is equivalent to $\varphi = 1$.
	A visualization of the results for $N = 16$, 32 and 64 is given in Fig. \ref{onewave_comp}.
	\begin{table}
		\caption{Relative $L^{2}$-error at $t=480$ for $2L$-periodic problem with initial data $u_{p1}$}
		\centering
		\begin{tabular}{ccc}
			\toprule
			$N$ & $E$ & rate \\
			\midrule
			16 & 0.14960222 & \multirow{2}{*}{2.41} \\
			32 & 0.02807195 & \multirow{2}{*}{2.28} \\
			64 & 0.00577740 & \multirow{2}{*}{2.16} \\
			128 & 0.00129088 & \multirow{2}{*}{2.07} \\
			256 & 0.00030683 &  \multirow{2}{*}{1.97} \\
			512 & 0.00007805 & \multirow{2}{*}{1.85} \\
			1024 & 0.00002172 & \\
		\end{tabular}
		\label{onewave_norms}
	\end{table}
	Again the plot indicates that the numerical approximation closes in on the exact solution and this is confirmed by the errors in Table \ref{onewave_norms} which are decreasing with a rate of approximately 2, as predicted by theory.
	The reason for this better behavior compared to the previous example could also be its somewhat less complicated nature, where the exact solution is simply the translation of a single solitary wave.
	
	\acknowledgement{The author is grateful to Helge Holden for his encouraging and support of the author's participation in the HYP2016 conference, and to Institut Mittag-Leffler for its hospitality during the Fall of 2016, providing an excellent working environment for this research.}
	
%
%
%

\end{document}